\documentclass[preprint,12pt]{elsarticle}

\usepackage{lineno}

\journal{arXiv (Vision 1)}

\usepackage{graphicx}
\usepackage{amssymb}
\usepackage{amstext}
\usepackage{amsmath}
\usepackage{amsthm}
\usepackage{bm}
\usepackage{float}
\usepackage{titlesec}
\usepackage{listings}
\usepackage{booktabs}
\usepackage{xcolor}
\usepackage{stfloats}
\usepackage{subfigure}
\usepackage{pgf}
\usepackage{wrapfig} 
\usepackage{tikz}
\usetikzlibrary{arrows}
\usepackage{appendix} 
\usepackage{algorithm}
\usepackage{algorithmic}
\usepackage{geometry}
\newtheorem{theorem}{Theorem}[section] 
\newtheorem{lemma}[theorem]{Lemma}
\theoremstyle{definition}
\newtheorem{definition}[theorem]{Definition}
\theoremstyle{definition}
\newtheorem{signmark}[theorem]{Mark}

\begin{document}
\begin{frontmatter}
\title{A General Verification for Functional Completeness by Abstract Operators}
\author[cor1]{Yang Tian\corref{cor1}}
\cortext[cor1]{Corresponding author}
\ead{tianyang16@mails.tsinghua.edu.cn}
\address{Tsinghua University, Beijing, China}
\begin{abstract}
An operator set is functionally incomplete if it can not represent the full set $\lbrace
\neg,\vee,\wedge,\rightarrow,\leftrightarrow\rbrace$. The verification for the functional incompleteness highly relies on constructive proofs. The judgement with a large untested operator set can be inefficient. Given with a mass of potential operators proposed in various logic systems, a general verification method for their functional completeness is demanded. This paper offers an universal verification for the functional completeness. Firstly, we propose two abstract operators $\widehat{R}$ and $\breve{R}$, both of which have no fixed form and are only defined by several weak constraints. Specially, $\widehat{R}_{\geq}$ and $\breve{R}_{\geq}$ are the abstract operators defined with the total order relation $\geq$. Then, we prove that any operator set $\mathfrak{R}$ is functionally complete if and only if it can represent the composite operator $\widehat{R}_{\geq}\circ\breve{R}_{\geq}$ or $\breve{R}_{\geq}\circ\widehat{R}_{\geq}$. Otherwise $\mathfrak{R}$ is determined to be functionally incomplete. This theory can be generally applied to any untested operator set to determine whether it is functionally complete.
\end{abstract}
\end{frontmatter}


\section{Introduction}
A set of truth function operators (or propositional connectives) is functionally complete if and only if all formulas constructed by $\lbrace
\neg,\vee,\wedge,\rightarrow,\leftrightarrow\rbrace$ can also be defined only based on that operator set. Otherwise the operator set is functionally incomplete \cite{shoenfield2018mathematical}. 

Determining whether a given operator set is functionally complete has both theoretical significance and practical value. For logic theories, verifying the functional completeness helps to analyse the mathematical properties of a given logic system defined with non-classical operators \cite{radzki2017axiom,belnap1970every,venema1990expressiveness}. For logic applications, the functionally complete set acts as the key unit in logic circuits and is the designed to realize demanded logic operation efficiently \cite{wang2016functionally,lehtonen2009stateful,varshavsky2004functionally,williamson1997resonant,yetter1993functionally,tokmen1978functionally}.

However, the definition of functional completeness can not provide a direct and efficient verification method for itself. For an untested operator set, the verification for its functional completeness usually relies on constructive proofs. There have been abundant researches that demonstrate object-specific verification proofs for the functional completeness of various truth function operator families. For a given truth function operator family, they develop self-contained and focused method to confirm the function completeness \cite{graham1967eta,massey1967binary,pinkava1978class,massey1970binary}. As for the much more general cases, a simple characterization method of functionally complete operator sets is proposed by researchers from \cite{maksimovic2006simple}, but the operator set is limited to be one-element.

To offer a general verification for the functional incompleteness of any given operator set, our work includes following parts:
\begin{itemize}
\item At first, we define two abstract truth function operators $\breve{R}$ and $\widehat{R}$ in \textbf{Subsection \ref{S21}}. They are abstract since they have no fixed forms and are defined only by several weak and elementary constraints. Thus, they can represent any operator satisfies those conditions. Specially, if the abstract operators are defined with the total relation $\geq$, we mark them as $\breve{R}_{\geq}$ and $\widehat{R}_{\geq}$.
\item At second, we analyse the properties of those two abstract operators $\breve{R}$ and $\widehat{R}$ in \textbf{Subsection \ref{S2}}. Then, in \textbf{Subsection \ref{S3}}, we show the connection between abstract operators and the known operators, where the relation between $\widehat{R}_{\geq}$ and $\neg$ is confirmed.
\item At third, we propose a new conception called semi-expressiveness in \textbf{Subsection \ref{S4}}. For any given operator set $\mathfrak{R}$, it's semi-expressive for $\lbrace\breve{R},\widehat{R}\rbrace$ if and only if it can represent $\widehat{R}\circ\breve{R}$ or $\breve{R}\circ\widehat{R}$. Moreover, we prove if $\mathfrak{R}$ is semi-expressive for $\lbrace\breve{R},\widehat{R}\rbrace$, then it can represent $\breve{R}$ and $\widehat{R}$.
\item Finally, we prove that any operator set $\mathfrak{R}$ is functionally complete if and only if it is semi-expressive for $\lbrace\breve{R}_{\geq},\widehat{R}_{\geq}\rbrace$ in \textbf{Subsection \ref{S5}}, which can offer a general verification for the functional completeness of any operator set. And in \textbf{Subsection \ref{S6}}, this theory is demonstrated to verify the functional completeness of several operator sets.
\end{itemize}
Our theory is demonstrated in a general but elementary form. All definitions and theorems are proposed to be self-contained.

\section{Review of elementary definitions}
To provide a clear vision, we propose a brief review for several elementary conceptions directly relevant with functional completeness. Later in the further introduction of our theory, all notations will follow the sign convention reviewed in this section.
\subsection{Truth function operators, truth value set and the valuation mapping}
Here we propose the sign convention for truth value set, truth function operators and the valuation mapping respectively.

The truth value set can be defined in various forms according to the demands of different logic systems. In our research, we only focus on the basic properties of it, which is given as following:
\begin{definition}
Let $V\subset [0,1]$ be the truth value set, which satisfies 
\begin{itemize}
\item $\lbrace 0,1\rbrace\subset V$;
\item $\forall v\in V$, there is $1-v \in V$.
\end{itemize} \label{ValueSet}
\end{definition}

Since we don't carry out object-specific analysis for any logic system, the truth function operator concerned in our research is proposed in a general form.
\begin{definition}
Let $\mathfrak{R}=\cup_{n}\mathfrak{R}^{n}$ be the set of all truth function operators needed to be analysed, where $\mathfrak{R}^{n}=\lbrace R\mid R:V^{n}\rightarrow V\rbrace$ denotes the set of all $n$-adic truth function operators in $\mathfrak{R}$.
\end{definition}

Given with a set of operators, we can represent formulas (or proportions) based on them. Thus, the synatx of a logic language can be defined.
\begin{definition}
Let $LP\left(\mathfrak{R}\right)$ be the set of all formulas constructed based on the truth function operators of $\mathfrak{R}$.
\end{definition}

The valuation mapping is used to assign truth value for a given formula and define the semantics for the logic language. From the perspective of algebra, it acts as the homomorphic mapping between the formula set to the truth value set.
\begin{definition}
Let $v:LP\left(\mathfrak{R}\right)\rightarrow V$ be the valuation function that determines the truth value of each formula in $LP\left(\mathfrak{R}\right)$. Note that $v$ is the homomorphic mapping between $LP\left(\mathfrak{R}\right)$ and $V$, such that for any $R$ of $\mathfrak{R}^{n}$, there is
\begin{equation}
\forall\phi_{1},\ldots,\phi_{n}\in LP\left(\mathfrak{R}\right),\; v\left(R\left(\phi_{1},\ldots,\phi_{n}\right)\right)=R\left(v\left(\phi_{1}\right),\ldots,v\left(\phi_{n}\right)\right).
\end{equation}
\end{definition}

\subsection{Functionally complete set}
After the synatx and semantics of the logic language is reviewed in general forms, we turn to defining the functional completeness.

\begin{definition}
Any operator set $\mathfrak{R}$ is functionally complete if and only if for every formula $\phi$ constructed using the full set of operators $\lbrace
\neg,\vee,\wedge,\rightarrow,\leftrightarrow\rbrace$, there is at least one formula $\psi$ of $LP\left(\mathfrak{R}\right)$ satisfies $\models\left(\phi\leftrightarrow \psi\right)$. \label{DefFCS}
\end{definition} 

For convenience, we refer the property proposed in \textbf{Definition \ref{DefFCS}} as that $\mathfrak{R}$ can represent $\lbrace
\neg,\vee,\wedge,\rightarrow,\leftrightarrow\rbrace$. The functional completeness of a truth function operator set is of great significance, since it determines how expressive a logic language defined based on that operator set is.

Moreover, in the studies relevant with functional completeness, the minimum functionally complete set often plays an important role, which is given by
\begin{definition}
For a given functionally complete set $\mathfrak{R}$, if any proper subset $\mathfrak{R}^{\prime}$ of it is functionally incomplete, then it is a minimum functionally complete set.
\end{definition}

In summary, we have reviewed the elementary conceptions related to the truth function operator set and the functional completeness. They will be frequently used in further analyses.

\section{Our work: a general verification for functional completeness}
We then turn to introducing our work. In our theory, we start with defining two new abstract $\breve{R}$ and $\widehat{R}$, then we propose the sufficient and necessary condition for a operator set $\mathfrak{R}$ to be functionally complete step by step.
\subsection{Abstract truth function operators $\breve{R}$ and $\widehat{R}$} \label{S21}
A general verification method should be able to fit any truth function operator set. Although we can not enumerate all possible truth function operators, we can still try to summarize their key properties closely related to the functional completeness. In our research, we use two abstract operators to contain those properties. An abstract operator is the truth function operator that has no fixed form and defined only by several elementary constraints, which can offer a general and abstract representation for a family of operators.

To make our work easy to follow, we start with giving brief definitions of those two abstract operators in this subsection. In later analyses, we will show the idea behind them and the relation between them and functional completeness systematically.

Before our analyses, for convenience, we define a sign mark related to the total order relation to simplify our introduction process.
\begin{signmark}
For any element $x$ and a non-empty set $Y$, if there exists a total order relation such that $x\triangleright y$ for each $y$ in $Y$, then we mark it as $x\triangleright Y$. Dually, we can define $Y\triangleright x$.
\end{signmark}

Then, there are two abstract truth function operators defined in our research, they are the choice operator $\breve{R}$ and the modification operator $\widehat{R}$.

\begin{definition}
We define the choice operator as $\breve{R}:V^{n}\rightarrow V$ with $n\geq 2$. It has no fixed form and is given by
\begin{equation}
v\left(\breve{R}\left(\phi_{1},\ldots,\phi_{n}\right)\right)=v\left(\phi_{i}\right),\;i\in\mathbb{Z}\cap\left[1,n\right],
\end{equation}
where there exists a total order relation $\triangleright$ such that $v\left(\phi_{i}\right)\triangleright\lbrace v\left(\phi_{j}\right)\mid j\neq i\rbrace$.\label{Choose}
\end{definition}
Based on \textbf{Definition \ref{Choose}}, it can be seen that the choice operator $\breve{R}$ is proposed in a very general form. We can simply treat $\lbrace v\left(\phi_{i}\right)\rbrace_{i\in I}$ as a set defined with a total order relation, then $\breve{R}$ select the greatest element in it.

\begin{definition}
We define the modification operator $\widehat{R}:V\rightarrow V$. It has no fixed form neither, and it satisfies if $v\left(\breve{R}\left(\phi_{1},\ldots,\phi_{n}\right)\right)=v\left(\phi_{i}\right)$ and $v\left(\phi_{i}\right)\triangleright\lbrace v\left(\phi_{j}\right)\mid j\neq i\rbrace$, then
\begin{itemize}
\item $\Big\lbrace v\left(\widehat{R}\left(\phi_{j}\right)\right)\mid j\neq i\Big\rbrace\triangleright v\left(\widehat{R}\left(\phi_{i}\right)\right)$;
\item $v\left(\left(\widehat{R}\circ\widehat{R}\right)\left(\phi_{i}\right)\right)=v\left(\phi_{i}\right)$;
\end{itemize} \label{Modify}
\end{definition}
Note that $\circ$ denotes the composite of functions. For example, $\left(\widehat{R}\circ\widehat{R}\right)\left(\phi_{i}\right)=\widehat{R}\left(\widehat{R}\left(\phi_{i}\right)\right)$. In \textbf{Definition \ref{Modify}}, the modification operator $\widehat{R}$ is defined dependently. A modification operator needs to correspond to at least one choice operator. Given that a choice operator can select the greatest element under a total order relation $\triangleright$, the modification operator can propose the inverse relation of $\triangleright$.

Apart of that, for convenience, we propose a sign convention for the modification operator $\widehat{R}$, which is given as
\begin{signmark}
When the modification operator $\widehat{R}$ is applied to a set of formulas $\lbrace\phi_{1},\ldots,\phi_{n}\rbrace$ respectively, we mark that $\widehat{R}\left(\phi_{1},\ldots,\phi_{n}\right)=\lbrace\widehat{R}\left(\phi_{1}\right),\ldots,\widehat{R}\left(\phi_{n}\right)\rbrace$.
\end{signmark}

\subsection{The properties of abstract truth function operators $\breve{R}$ and $\widehat{R}$} \label{S2}

After defining abstract truth function operators $\breve{R}$ and $\widehat{R}$, here we analyse several basic characters of them.

Firstly, since the modification operator $\widehat{R}$ defined in \textbf{Definition \ref{Modify}} can propose the inverse relation for a given total order relation $\triangleright$, it's easy to know the following theorem.
\begin{theorem}
If there exists a total order relation $\triangleright$ such that $v\left(\phi\right)\triangleright v\left(\psi\right)$, then $v\left(\widehat{R}\left(\psi\right)\right)\triangleright v\left(\widehat{R}\left(\phi\right)\right)$. \label{Theorem0}
\end{theorem}
\begin{proof}
For the choice operator $\breve{R}:V^{n}\rightarrow V$, assume that $n=2$. Based on \textbf{Definition \ref{Choose}}, if $v\left(\phi\right)\triangleright v\left(\psi\right)$, then $v\left(\breve{R}\left(\phi,\psi\right)\right)=v\left(\phi\right)$. Based on \textbf{Definition \ref{Modify}}, there is $v\left(\widehat{R}\left(\psi\right)\right)\triangleright v\left(\widehat{R}\left(\phi\right)\right)$.
\end{proof}

Secondly, being able to define the inverse relation for a given total order relation, the modification operator $\widehat{R}$ itself can be proved as a bijection.
\begin{theorem}
The modification operator $\widehat{R}:V\rightarrow V$ is a bijection from $V$ to $V$. \label{Theorem2}
\end{theorem}
\begin{proof}
Our proof is demonstrated with two steps.
\begin{itemize}
\item At first, we prove that $\widehat{R}$ is injective by reduction to absurdity. Assume $\widehat{R}$ is not injective, then
\begin{equation}
\exists v_{i}\neq v_{j}\in V,\;\widehat{R}\left(v_{i}\right)=\widehat{R}\left(v_{j}\right)=v_{k}, \label{9}
\end{equation}
based on \textbf{Definition \ref{Modify}}, there is
\begin{align}
v_{i}=\left(\widehat{R}\circ\widehat{R}\right)\left(v_{i}\right)=\widehat{R}\left(v_{k}\right)=\left(\widehat{R}\circ\widehat{R}\right)\left(v_{j}\right)=v_{j}, \label{10}
\end{align}
which is contradicts with \textbf{(\ref{9})}. Thus, $\widehat{R}$ must be injective.
\item At second, we prove that $\widehat{R}$ is a surjection. For convenience, we mark the domain of definition as $\mathtt{dom}\;\widehat{R}$ and mark the value range as $\mathtt{ran}\;\widehat{R}$. Assume that $\widehat{R}$ is not a surjection, then there exist $v_{j}\in \mathtt{ran}\;\widehat{R}$ such that for all $v_{i}\in \mathtt{dom}\;\widehat{R}$, $\widehat{R}\left(v_{i}\right)\neq v_{j}$. Then, since $ \mathtt{dom}\;\widehat{R}= \mathtt{ran}\;\widehat{R}$, there must be 
\begin{equation}
\exists v_{k}\in \mathtt{ran}\;\widehat{R},\; \vert\lbrace v_{i}\mid \widehat{R}\left(v_{i}\right)= v_{k}\rbrace\vert>1, \label{11}
\end{equation}
based on \textbf{Definition \ref{Modify}}, \textbf{(\ref{11})} can implies the same contradiction in \textbf{(\ref{10})}. Thus, $\widehat{R}$ must be a surjection.
\end{itemize}
Combine those two conclusions, we can finish our proof. 
\end{proof}

\textbf{Theorem \ref{Theorem0}} shows an important property of the modification operator $\widehat{R}$. Given that it's a bijection for $V$ to $V$, there are a lot of theorems can be applied to it to simplify our proofs (such as the proof of \textbf{Theorem \ref{Theorem4}}).

Thirdly, we also pay special attention to the properties of the composite operators $\breve{R}\circ\widehat{R}$ and $\widehat{R}\circ\breve{R}$, in the following theorem, we will analyse what those composite operators can imply on the truth value of given formulas.
\begin{theorem}
If $v\left(\breve{R}\left(\phi_{1},\ldots,\phi_{n}\right)\right)=v\left(\phi_{i}\right)$ and $v\left(\phi_{i}\right)\triangleright\lbrace v\left(\phi_{j}\right)\mid j\neq i\rbrace$, then we know
\begin{itemize}
\item assume $v\left(\left(\breve{R}\circ\widehat{R}\right)\left(\phi_{1},\ldots,\phi_{n}\right)\right)=v\left(\widehat{R}\left(\phi_{k}\right)\right)$, then $v\left(\phi_{i}\right)\triangleright\lbrace v\left(\phi_{j}\right)\mid j\neq i\;\text{and}\;j\neq k\rbrace\triangleright v\left(\phi_{k}\right)$;
\item there is $v\left(\left(\widehat{R}\circ\breve{R}\right)\left(\phi_{1},\ldots,\phi_{n}\right)\right)=v\left(\widehat{R}\left(\phi_{i}\right)\right)$.
\end{itemize}\label{Theorem1}
\end{theorem}

\begin{proof}
We respectively prove those two properties.
\begin{itemize}
\item For the first property, it's clear that 
\begin{align}
v\left(\left(\breve{R}\circ\widehat{R}\right)\left(\phi_{1},\ldots,\phi_{n}\right)\right)&=v\left(\breve{R}\left(\widehat{R}\left(\phi_{1}\right),\ldots,\widehat{R}\left(\phi_{n}\right)\right)\right),\\
&=v\left(\widehat{R}\left(\phi_{k}\right)\right),
\end{align}
based on \textbf{Definition \ref{Choose}}, it's clear that $v\left(\widehat{R}\left(\phi_{k}\right)\right)\triangleright\Big\lbrace v\left(\widehat{R}\left(\phi_{j}\right)\right)\mid j\neq k\Big\rbrace$. Given \textbf{Theorem \ref{Theorem1}}, we can know
\begin{equation}
\forall j\neq k,\;v\left(\left(\widehat{R}\circ\widehat{R}\right)\left(\phi_{j}\right)\right)\triangleright v\left(\left(\widehat{R}\circ\widehat{R}\right)\left(\phi_{k}\right)\right),
\end{equation}
which implies that 
\begin{align}
\Big\lbrace v\left(\left(\widehat{R}\circ\widehat{R}\right)\left(\phi_{j}\right)\right)\mid j\neq k\Big\rbrace &\triangleright v\left(\left(\widehat{R}\circ\widehat{R}\right)\left(\phi_{k}\right)\right),\\
\lbrace v\left(\phi_{j}\right)\mid j\neq k\rbrace &\triangleright v\left(\phi_{k}\right).
\end{align}
Since $v\left(\phi_{i}\right)\triangleright\lbrace v\left(\phi_{j}\right)\mid j\neq i\rbrace$, there is
\begin{equation}
v\left(\phi_{i}\right)\triangleright\lbrace v\left(\phi_{j}\right)\mid j\neq i\;\text{and}\;j\neq k\rbrace\triangleright v\left(\phi_{k}\right),
\end{equation}
which finishes the proof.
\item The proof of the second property can be directly obtained based on that $v\left(\breve{R}\left(\phi_{1},\ldots,\phi_{n}\right)\right)=v\left(\phi_{i}\right)$.
\end{itemize}
\end{proof}

\subsection{The connection between abstract operators and known operators} \label{S3}
 In \textbf{Subsection \ref{S21}}, we have suggested that $\breve{R}$ and $\widehat{R}$ can be treated as the general and abstract representation for a family of operators with fixed forms. Up to now, we have analysed several elementary properties of those two abstract operators. We then turn to indicating the underlying connection between abstract operators and other known operators.

Those two abstract operators are defined based on the total order relation. Since we need to talk about different known operators defined in fixed forms, we define a sign convention to distinguish between different total order relations.
\begin{signmark}
To distinguish between different total order relations, we define if there exists a total order relation $\triangleright$ such that $v\left(\breve{R}\left(\phi_{1},\ldots,\phi_{n}\right)\right)=v\left(\phi_{i}\right),\;i\in\mathbb{Z}\cap\left[1,n\right]$ and $v\left(\phi_{i}\right)\triangleright\lbrace v\left(\phi_{j}\right)\mid j\neq i\rbrace$, then the choice operator is marked as $\breve{R}_{\triangleright}$. And the corresponding modification operator is marked as $\widehat{R}_{\triangleright}$.
\end{signmark}

We are specially interested in the situation where those two abstract operators are defined with $\geq$ (greater than or equal to) and $\leq$ (less than or equal to), both of which are the most basic total order relations in the set of real numbers.
\begin{theorem}
If the choice operator is $\breve{R}_{\geq}$, then the corresponding modify operator $\widehat{R}_{\geq}$ must satisfy 
\begin{equation}
\widehat{R}_{\geq}\left(v\right)=1-v. \label{KCO}
\end{equation}
And the same conclusion can be obtained when the choice operator is $\breve{R}_{\leq}$.\label{Theorem4}
\end{theorem}

\begin{proof}
To prove \textbf{Theorem \ref{Theorem4}}, we can prove that \textbf{(\ref{KCO})} is the sufficient and necessary condition for that $\widehat{R}_{\geq}$ satisfies the \textbf{Theorem \ref{Theorem0}}.
\begin{itemize}
\item At first, the sufficiency is trivial to prove based on
\begin{equation}
\forall v_{i}\geq v_{j}\in V,\;\left(1-v_{j}\right)\geq \left(1-v_{i}\right).\label{KCOS}
\end{equation}
\item At second, we prove the necessity by reduction to absurdity. We assume a case where 
\begin{align}
\exists v_{i}\in V,\;k\neq 0,\;\widehat{R}_{\geq}\left(v_{i}\right)=1-v_{i}+k\in V, \label{Assume1}
\end{align}
which implies that $\widehat{R}_{\geq}\left(v_{i}\right)$ doesn't satisfy \textbf{(\ref{KCO})}. And based on \textbf{Definition \ref{ValueSet}}, it can be known that $v_{i}-k\in V$. Then, we further assume that 
\begin{equation}
\forall v_{j}\neq v_{i} \neq v_{i}-k\in V,\; \widehat{R}_{\geq}\left(v_{j}\right)=1-v_{j}. \label{Assume2}
\end{equation}
Based on \textbf{(\ref{Assume1})} and \textbf{(\ref{Assume2})}, we can define a situation where $\widehat{R}_{\geq}$ can satisfy \textbf{(\ref{KCO})} when acts on all elements of $V$, except for $v_{i}$ and $v_{i}-k$. Since we have proved that $\widehat{R}_{\geq}$ is a bijection from $V$ to $V$ in \textbf{Theorem \ref{Theorem2}}, given with \textbf{(\ref{Assume1})} and \textbf{(\ref{Assume2})}, there must be
\begin{equation}
\widehat{R}_{\geq}\left(v_{i}-k\right)=1-v_{i}.
\end{equation}
And it's easy to know
\begin{equation}
\begin{cases}
\widehat{R}_{\geq}\left(v_{i}\right)\geq\widehat{R}_{\geq}\left(v_{i}-k\right), & v_{i}\geq v_{i}-k\\
\widehat{R}_{\geq}\left(v_{i}-k\right)\geq\widehat{R}_{\geq}\left(v_{i}\right), & v_{i}-k\geq v_{i}
\end{cases},
\end{equation}
which implies the contradiction with \textbf{Theorem \ref{Theorem0}}. Thus, the case assumed in \textbf{(\ref{Assume1})} can not happen.
\end{itemize}
To sum up, we know that $\widehat{R}_{\geq}\left(v\right)=1-v$ is the sufficient and necessary condition for that $\widehat{R}_{\geq}$ satisfies \textbf{Theorem \ref{Theorem0}}, which finishes the proof. Following those proposed steps, we can also prove the conclusion for $\widehat{R}_{\leq}$.
\end{proof}

\textbf{Theorem \ref{Theorem4}} shows a simple and clear result: when abstract operators are defined with $\geq$ and $\leq$, the modification operator has same function as the classical operator $\neg$ (we know $\neg\left(v\right)=1-v$). This inspires us to think about whether those two abstract operators can be used to $\lbrace
\neg,\vee,\wedge,\rightarrow,\leftrightarrow\rbrace$ when they are defined with specific total order relations.

\subsection{The semi-expressiveness for abstract operators} \label{S4}
Before analysing the relation between the abstract operators and the functional completeness, here we propose a new conception that is called semi-expressiveness.

\begin{definition}
A operator set $\mathfrak{R}$ is semi-expressive for $\lbrace\widehat{R},\breve{R}\rbrace$ if there exist at least one set of operators $R_{1},\ldots,R_{k}$ in it satisfy \label{SEMIE}
\begin{equation}
v\Big(\left(R_{1}\circ\ldots\circ R_{k}\right)\left(\phi_{1},\ldots,\phi_{n}\right)\Big)=v\left(\left(\widehat{R}\circ \breve{R}\right)\left(\phi_{1},\ldots,\phi_{n}\right)\right), \label{Con1}
\end{equation}
or
\begin{equation}
v\Big(\left(R_{1}\circ\ldots\circ R_{k}\right)\left(\phi_{1},\ldots,\phi_{n}\right)\Big)=v\left(\left(\breve{R}\circ\widehat{R}\right)\left(\phi_{1},\ldots,\phi_{n}\right)\right). \label{Con2}
\end{equation}
\end{definition}

In brief, a operator set $\mathfrak{R}$ is semi-expressive for $\lbrace\widehat{R},\breve{R}\rbrace$ means that merely based on the operators $\mathfrak{R}$ contains, the composite operators $\widehat{R}\circ \breve{R}$ and $\breve{R}\circ\widehat{R}$ can be represented. 

The choice operator $\breve{R}$ can select the greatest element in a given subset of $V$ with a given total order relation $\triangleright$, and the modification operator can propose the inverse relation $\triangleleft$ of $\triangleright$. So we can see that 
\begin{itemize}
\item if \textbf{(\ref{Con1})} is satisfied, then greatest element of a given subset of $V$ is selected under the total order relation $\triangleright$ and modified to be the least element;
\item if \textbf{(\ref{Con2})} is satisfied, then least element of a given subset of $V$ is selected under the total order relation $\triangleright$ (the greatest element under $\triangleleft$ is the least element under $\triangleright$) and modified to be the greatest element.
\end{itemize}
In a word, a operator set $\mathfrak{R}$ that is semi-expressive for $\lbrace\widehat{R},\breve{R}\rbrace$ can realize the transformation between the least and greatest element in a given subset of $V$ with a total order relation.

Just as what have been mentioned for $\lbrace\widehat{R},\breve{R}\rbrace$, the semi-expressiveness is the ability to represent the composite operators $\widehat{R}\circ \breve{R}$ and $\breve{R}\circ\widehat{R}$. It's natural to wander if it also means the ability to represent $\widehat{R}$ and $\breve{R}$ respectively.

\begin{theorem}
If $\mathfrak{R}$ is semi-expressive for $\lbrace\widehat{R},\breve{R}\rbrace$, then $\mathfrak{R}$ can represent $\widehat{R}$ and $\breve{R}$. \label{Respectively}
\end{theorem}

\begin{proof}
The following is our proof.
\begin{itemize}
\item For $\breve{R}$,
\begin{itemize}
\item if \textbf{(\ref{Con1})} is satisfied, we assume that $\breve{R}\left(\phi_{1},\ldots,\phi_{n}\right)=\phi_{i}$, then
\begin{align}
&\left(\widehat{R}\circ \breve{R}\right)\left(\left(\widehat{R}\circ \breve{R}\right)\left(\phi_{1},\ldots,\phi_{n}\right),\ldots,\left(\widehat{R}\circ \breve{R}\right)\left(\phi_{1},\ldots,\phi_{n}\right)\right)\\
=&\left(\widehat{R}\circ \breve{R}\right)\left(\widehat{R}\left(\phi_{i}\right),\ldots,\widehat{R}\left(\phi_{i}\right)\right)\\
=&\left(\widehat{R}\circ \widehat{R}\right)\left(\phi_{i}\right)\\
=&\breve{R}\left(\phi_{1},\ldots,\phi_{n}\right). \label{23}
\end{align}
Note that \textbf{(\ref{23})} is implied by that for any formula $\phi_{i}$, $\left(\widehat{R}\circ \widehat{R}\right)\left(\phi_{i}\right)=\phi_{i}$, which is defined in \textbf{Definition \ref{Modify}}.
\item Else if \textbf{(\ref{Con2})} is satisfied, there is
\begin{align}
&\left(\breve{R}\circ\widehat{R}\right)\left(\left(\breve{R}\circ\widehat{R}\right)\left(\phi_{1},\ldots,\phi_{1}\right),\ldots,\left(\breve{R}\circ\widehat{R}\right)\left(\phi_{n},\ldots,\phi_{n}\right)\right)\\
=&\left(\breve{R}\circ\widehat{R}\right)\left(\widehat{R}\left(\phi_{1}\right),\ldots,\widehat{R}\left(\phi_{n}\right)\right)\\
=&\left[\breve{R}\circ\left(\widehat{R}\circ \widehat{R}\right)\right]\left(\phi_{1},\ldots,\phi_{n}\right)\\
=&\breve{R}\left(\phi_{1},\ldots,\phi_{n}\right).
\end{align}
\end{itemize}
Based on what have been proved, we know if $\mathfrak{R}$ is semi-expressive for $\lbrace\widehat{R},\breve{R}\rbrace$, then $\mathfrak{R}$ can represent $\breve{R}$.
\item For $\widehat{R}$, we know that no matter \textbf{(\ref{Con1})} or \textbf{(\ref{Con2})} is satisfied, there is
\begin{equation}
\forall \phi_{i},\;\left(\widehat{R}\circ \breve{R}\right)\left(\phi_{i},\ldots\phi_{i}\right)=\left(\breve{R}\circ\widehat{R}\right)\left(\phi_{i},\ldots\phi_{i}\right)=\widehat{R}\left(\phi_{i}\right)
\end{equation} 
Thus, if $\mathfrak{R}$ is semi-expressive for $\lbrace\widehat{R},\breve{R}\rbrace$, then $\mathfrak{R}$ can represent $\widehat{R}$.
\end{itemize}
To summarize, we can finally draw a conclusion that \textbf{Theorem \ref{SEMIE}} is true.
\end{proof}

Up to now, we have prepared all the necessary definitions and theorems about the abstract operators. Later in the next subsection, we will use them to propose a sufficient and necessary condition of functional completeness.

\subsection{A precise verification for the functional completeness} \label{S5}

Here we propose our theorem of the sufficient and necessary condition of functional completeness.
\begin{theorem}
For a operator set $\mathfrak{R}$, it is functionally complete if and only if it is semi-expressive for $\lbrace\widehat{R}_{\geq},\breve{R}_{\geq}\rbrace$. \label{Theorem3}
\end{theorem}

\begin{proof}
To prove the sufficient and necessary condition proposed in \textbf{Theorem \ref{Theorem3}}, our proof consists of two parts. 
\begin{itemize}
\item At first, based on \textbf{Theorem \ref{SEMIE}}, we know that if $\mathfrak{R}$ is semi-expressive for $\lbrace\widehat{R}_{\geq},\breve{R}_{\geq}\rbrace$, then $\mathfrak{R}$ can represent $\widehat{R}_{\geq}$ and $\breve{R}_{\geq}$. 
\item At second, we prove the sufficiency, which means that every operator set that is semi-expressive for $\lbrace\widehat{R}_{\geq},\breve{R}_{\geq}\rbrace$ will be functionally complete. Given that $\mathfrak{R}$ can represent $\widehat{R}_{\geq}$ and $\breve{R}_{\geq}$, we can have following analyses.
\begin{itemize}
\item For $\breve{R}_{\geq}$, we know if $v\left(\breve{R}_{\geq}\left(\phi_{1},\ldots,\phi_{n}\right)\right)=v\left(\phi_{i}\right)$, and $v\left(\phi_{i}\right)\geq\lbrace v\left(\phi_{j}\right)\mid j\neq i\rbrace$, then it's trivial that $v\left(\phi_{i}\right)=v\Big(\left(\left(\left(\phi_{1}\vee\phi_{2}\right)\ldots\right)\vee\phi_{n}\right)\Big)$. Thus, if $\mathfrak{R}$ represents $\breve{R}_{\geq}$, then $\mathfrak{R}$ can represent $\vee$. 
\item For $\widehat{R}_{\geq}$, since for each formula $\phi_{i}$ $v\left(\neg \phi_{i}\right)=1-v\left(\phi_{i}\right)$, we can prove if $\mathfrak{R}$ represents $\widehat{R}_{\geq}$, then $\mathfrak{R}$ can represent $\neg$ directly based on \textbf{Theorem \ref{Theorem4}}.
\end{itemize}
Based on what have be proved, we know a operator set $\mathfrak{R}$ that is semi-expressive for $\lbrace\widehat{R}_{\geq},\breve{R}_{\geq}\rbrace$ can represent $\lbrace \neg,\vee\rbrace$, which is a minimum functionally complete operator set. It's known if $\mathfrak{R}$ can represent a functionally complete operator set, then it can represent $\lbrace\neg,\vee,\wedge,\rightarrow,\leftrightarrow\rbrace$. Thus, the sufficiency has been proved.

\item At third, we prove the necessity, which means that every functionally complete operator set must be semi-expressive for $\lbrace\widehat{R}_{\geq},\breve{R}_{\geq}\rbrace$. Following the idea used in \textbf{Theorem \ref{SEMIE}}, we know that 
\begin{equation}
\forall \phi_{i},\;\left(\widehat{R}_{\geq}\circ \breve{R}_{\geq}\right)\left(\phi_{i},\ldots\phi_{i}\right)=\left(\breve{R}_{\geq}\circ\widehat{R}_{\geq}\right)\left(\phi_{i},\ldots\phi_{i}\right)=\widehat{R}_{\geq}\left(\phi_{i}\right)
\end{equation} 
if there doesn't exist any set of operators $R_{1},\ldots,R_{k}$ in $\mathfrak{R}$ satisfies \textbf{(\ref{Con1})} nor \textbf{(\ref{Con2})}, then for any formula $\phi$, $v\left(\widehat{R}_{\geq}\left(\phi\right)\right)$ can never be defined by $\mathfrak{R}$. Based on \textbf{Theorem \ref{Theorem4}}, it's easy to prove that $\neg$ can't be represented by $\mathfrak{R}$. Thus, the necessity is proved.
\end{itemize}

To sum up, we have proved that a operator set $\mathfrak{R}$ is functionally complete if and only if it is semi-expressive for $\lbrace\widehat{R}_{\geq},\breve{R}_{\geq}\rbrace$.
\end{proof}

Note that since $\geq$ and $\leq$ are inverse relations of each other, we can replace $\lbrace\widehat{R}_{\geq},\breve{R}_{\geq}\rbrace$ in \textbf{Theorem \ref{Theorem3}} as $\lbrace\widehat{R}_{\leq},\breve{R}_{\leq}\rbrace$ and it still holds.

The \textbf{Theorem \ref{Theorem3}} and its proof are simple since we have proposed all necessary foundations of them in the previous analysis. An important thing about \textbf{Theorem \ref{Theorem3}} is that the method can verify the functional completeness of any $n$-element untested operator set since there is no further limitation for it.

\subsection{Several demonstrations of our theory} \label{S6}

Finally, we propose two applications of \textbf{Theorem \ref{Theorem3}}, where we use it to verify the functional completeness of several truth function operator sets.

To provide a better understanding for the functional completeness verification based on \textbf{Theorem \ref{Theorem3}}, we propose several trivial examples at first.

\begin{lemma}
Define four operators \label{ASL}
\begin{align}
R_{1}\left(\phi_{1},\ldots,\phi_{n}\right)&=\neg\left(\left(\left(\phi_{1}\vee\phi_{2}\right)\ldots\right)\vee\phi_{n}\right),\label{30}\\
R_{2}\left(\phi_{1},\ldots,\phi_{n}\right)&=\neg\left(\left(\left(\phi_{1}\wedge\phi_{2}\right)\ldots\right)\wedge\phi_{n}\right),\label{31}\\
R_{3}\left(\phi_{1},\ldots,\phi_{n}\right)&=\left(\left(\left(\neg\phi_{1}\right)\vee\left(\neg\phi_{2}\right)\right)\ldots\right)\vee\left(\neg\phi_{n}\right),\label{32}\\
R_{4}\left(\phi_{1},\ldots,\phi_{n}\right)&=\left(\left(\left(\neg\phi_{1}\right)\wedge\left(\neg\phi_{2}\right)\right)\ldots\right)\wedge\left(\neg\phi_{n}\right),\label{33}
\end{align}
then $\lbrace R_{1}\rbrace$, $\lbrace R_{2}\rbrace$, $\lbrace R_{3}\rbrace$ and $\lbrace R_{4}\rbrace$ are all functionally complete.
\end{lemma}
\begin{proof}
It's trivial that \textbf{(\ref{30})} and \textbf{(\ref{31})} have the same form with \textbf{(\ref{Con1})} while \textbf{(\ref{32})} and \textbf{(\ref{33})} have the same form with \textbf{(\ref{Con2})}, so \textbf{Lemma \ref{ASL}} can be directly proved based on \textbf{Theorem \ref{Theorem3}}.
\end{proof}

Then, we propose three elementary but non-trivial examples.

\begin{lemma}
Define three operators \label{ASL2}
\begin{align}
R_{5}\left(\phi_{1},\ldots,\phi_{n}\right)&=\neg\left(\left(\left(\phi_{1}\bigstar_{1}\phi_{2}\right)\ldots\right)\bigstar_{n-1}\phi_{n}\right),\label{34}\\
R_{6}\left(\phi_{1},\ldots,\phi_{n}\right)&=\left(\left(\left(\neg\phi_{1}\right)\bigstar_{1}\left(\neg\phi_{2}\right)\right)\ldots\right)\bigstar_{n-1}\left(\neg\phi_{n}\right),\label{35}\\
R_{7}\left(\phi_{1},\ldots,\phi_{n}\right)&=\left(\left(\left(\blacklozenge_{1}\phi_{1}\right)\bigstar_{1}\left(\blacklozenge_{2}\phi_{2}\right)\right)\ldots\right)\bigstar_{n-1}\left(\blacklozenge_{n}\phi_{n}\right),\label{36}
\end{align}
where each $\bigstar_{k}$ is randomly selected from $\lbrace\vee,\wedge\rbrace$ and each $\blacklozenge_{k}$ is randomly selected from $\lbrace\neg,\neg\neg\rbrace$. Then 
\begin{itemize}
\item $\lbrace R_{5}\rbrace$ and $\lbrace R_{6}\rbrace$ are all functionally complete;
\item $\lbrace R_{7}\rbrace$ is functionally incomplete if there exists at least one $x$ ($x\leq n$) that satisfies $\blacklozenge_{x}$ is $\neg\neg$.
\end{itemize}
\end{lemma}

\begin{proof}
We respectively prove \textbf{(\ref{34})}, \textbf{(\ref{35})} and \textbf{(\ref{36})}.
\begin{itemize}
\item For \textbf{(\ref{34})}, given that each $\bigstar_{k}$ is randomly selected from $\lbrace\vee,\wedge\rbrace$, our idea to demonstrate the proof is given as following:
\begin{itemize}
\item with the truth values of given formulas $\lbrace v\left(\phi_{1}\right),\ldots,v\left(\phi_{n}\right)\rbrace$, we can find the greatest element $v\left(\phi_{i}\right)$ and least element $v\left(\phi_{j}\right)$ of this set with the total order relation $\geq$. 
\item Then for a given operator $R_{5}$, we can confirm whether $\bigstar_{n-1}$ is $\vee$ or $\wedge$, and
\begin{itemize}
\item if $\bigstar_{n-1}$ is $\vee$, then we know
\begin{equation}
v\left(\left(\widehat{R}_{\geq}\circ \breve{R}_{\geq}\right)\left(\phi_{1},\ldots,\phi_{n}\right)\right)=v\Big(\left(\neg\left(\left(\left(\ldots\bigstar_{1}\ldots\right)\ldots\right)\bigstar_{n-1}\phi_{i}\right)\right)\Big), \label{37}
\end{equation}
where $\phi_{i}$ is put at the end of the sequence so that $\bigstar_{n-1}$ can act on it. As for other formulas $\lbrace \phi_{x}\mid x\neq i\rbrace$, there is no limitation for their location so they are represented by ``$\ldots$" (this omission is also used in the later analysis). We can see that in \textbf{(\ref{37})}, $\lbrace R_{5}\rbrace$ is semi-expressive for $\lbrace\widehat{R}_{\geq},\breve{R}_{\geq}\rbrace$. Based on \textbf{Theorem \ref{Theorem3}}, it's functionally complete.
\item if $\bigstar_{n-1}$ is $\wedge$, then we know
\begin{equation}
v\left(\left(\widehat{R}_{\leq}\circ \breve{R}_{\leq}\right)\left(\phi_{1},\ldots,\phi_{n}\right)\right)=v\Big(\left(\neg\left(\left(\left(\ldots\bigstar_{1}\ldots\right)\ldots\right)\bigstar_{n-1}\phi_{j}\right)\right)\Big), \label{38}
\end{equation}
thus $\lbrace R_{5}\rbrace$ is semi-expressive for $\lbrace\widehat{R}_{\leq},\breve{R}_{\leq}\rbrace$, which makes $\lbrace R_{5}\rbrace$ be functionally complete.
\end{itemize}
\item To sum up, using the idea described above, we can conclude that $\lbrace R_{5}\rbrace$ described in \textbf{(\ref{34})} is functionally complete.
\end{itemize} 
\item For \textbf{(\ref{35})}, the same idea used to proof \textbf{(\ref{34})} can also be used. It's trivial that
\begin{itemize}
\item if $\bigstar_{n-1}$ is $\vee$, then 
\begin{equation}
v\left(\left(\breve{R}_{\geq}\circ \widehat{R}_{\geq}\right)\left(\phi_{1},\ldots,\phi_{n}\right)\right)=v\Big(\left(\left(\ldots\bigstar_{1}\ldots\right)\ldots\right)\bigstar_{n-1}\left(\neg\phi_{j}\right)\Big), \label{39}
\end{equation}
which implies that $\lbrace R_{6}\rbrace$ is functionally complete when it's semi-expressive for $\lbrace\widehat{R}_{\geq},\breve{R}_{\geq}\rbrace$.
\item if $\bigstar_{n-1}$ is $\wedge$, then 
\begin{equation}
v\left(\left(\breve{R}_{\leq}\circ \widehat{R}_{\leq}\right)\left(\phi_{1},\ldots,\phi_{n}\right)\right)=v\Big(\left(\left(\ldots\bigstar_{1}\ldots\right)\ldots\right)\bigstar_{n-1}\left(\neg\phi_{i}\right)\Big), \label{40}
\end{equation}
which implies that $\lbrace R_{6}\rbrace$ is semi-expressive for $\lbrace\widehat{R}_{\leq},\breve{R}_{\leq}\rbrace$ and functionally complete.
\end{itemize}
Thus, $\lbrace R_{6}\rbrace$ described in \textbf{(\ref{35})} is functionally complete.
\item For \textbf{(\ref{36})}, we know 
\begin{align}
v\Big(R_{7}\left(\phi_{1},\ldots,\phi_{n}\right)\Big)=v\left(\left(\breve{R}_{\geq}\circ \widehat{R}_{\geq}\right)\left(\phi_{1},\ldots,\phi_{n}\right)\right)\Leftrightarrow\\
 v\Big(\left(\left(\left(\blacklozenge_{1}\ldots\right)\bigstar_{1}\left(\blacklozenge_{2}\ldots\right)\right)\ldots\right)\bigstar_{n-1}\left(\blacklozenge_{n}\phi_{j}\right)\Big)=v\left(\neg\phi_{j}\right),\label{42}
\end{align}
or 
\begin{align}
v\Big(R_{7}\left(\phi_{1},\ldots,\phi_{n}\right)\Big)=v\left(\left(\breve{R}_{\leq}\circ \widehat{R}_{\leq}\right)\left(\phi_{1},\ldots,\phi_{n}\right)\right)\Leftrightarrow\\
 v\Big(\left(\left(\left(\blacklozenge_{1}\ldots\right)\bigstar_{1}\left(\blacklozenge_{2}\ldots\right)\right)\ldots\right)\bigstar_{n-1}\left(\blacklozenge_{n}\phi_{i}\right)\Big)=v\left(\neg\phi_{i}\right), \label{43}
\end{align}
However, if there exists at least one $x$ ($x\leq n$) that satisfies $\blacklozenge_{x}$ is $\neg\neg$, then 
\begin{itemize}
\item if every $\blacklozenge_{x}$ in \textbf{(\ref{36})} is $\neg\neg$, then for each formula $\phi_{x}$, there is $\blacklozenge_{x}\phi_{x}=\phi_{x}$. Thus, it's easy to know $\lbrace R_{7}\rbrace$ can never represent $\widehat{R}_{\leq}$ nor $\widehat{R}_{\geq}$. In other words, $\neg$ can not be represented by $\lbrace R_{7}\rbrace$, which makes $\lbrace R_{7}\rbrace$ be functionally incomplete;
\item if $\blacklozenge_{n}$ in \textbf{(\ref{36})} is $\neg$, then 
\begin{itemize}
\item assume $\bigstar_{n-1}$ is $\vee$, given that $v\left(\neg\phi_{j}\right)\geq\lbrace v\left(\neg\phi_{k}\right)\mid k\neq j\rbrace$, then it's trivial that \textbf{(\ref{42})} is guaranteed (no matter how each $\bigstar_{k}$ is selected) if and only if $v\left(\neg\phi_{j}\right)\geq\lbrace v\left(\phi_{k}\right)\mid k\neq j\rbrace$. For the formula set $\lbrace\phi_{1},\ldots,\phi_{n}\rbrace$ does not meet this condition, \textbf{(\ref{42})} can be false. Thus, $\lbrace R_{7}\rbrace$ is functionally incomplete;
\item assume $\bigstar_{n-1}$ is $\wedge$, it's trivial that \textbf{(\ref{43})} is guaranteed (no matter how each $\bigstar_{k}$ is selected) if and only if $v\left(\neg\phi_{i}\right)\leq\lbrace v\left(\phi_{k}\right)\mid k\neq i\rbrace$. For other situations, \textbf{(\ref{43})} can be false, which makes $\lbrace R_{7}\rbrace$ be functionally incomplete.
\end{itemize}
\end{itemize}
Finally, we can conclude that $\lbrace R_{7}\rbrace$ described in \textbf{(\ref{36})} is functionally incomplete with the given condition.
\end{itemize}
\end{proof}

The examples used in \textbf{Lemma \ref{ASL2}} are proposed in general forms. In fact, the well known Peirce operator $\downarrow$ (there is $\phi\downarrow\psi=\neg\left(\phi\vee\psi\right)$, $\lbrace\downarrow\rbrace$ is functionally complete \cite{brady2000peirce}) is a special case of \textbf{(\ref{34})}. For those non-trivial examples, the method proposed by us can provide an unified process to verify the functional completeness. 

\section{Conclusion and Discussion}
Let's have a brief review for what have been proposed in our theory. Our work aims at providing a general verification method for the functional completeness (or functional incompleteness). To realize this ambition, we define two abstract operators $\widehat{R}$ and $\breve{R}$ in \textbf{Definition \ref{Choose}} and \textbf{Definition \ref{Modify}} to construct a new theory, where we propose the sufficient and necessary condition for a truth function operator set to be functionally complete in \textbf{Theorem \ref{Theorem3}}. In our work, we demonstrate the theory on several examples, which suggests that our method can offer an efficient verification. Compared with the object-specific verification in previous \cite{graham1967eta,massey1967binary,pinkava1978class,massey1970binary}, our theory can be applied more generally.

In our theory, what interests us most is the properties of those two abstract operators proposed in \textbf{Subsection \ref{S2}}. Although we mainly use them to study the functional completeness, they might also be inspiring for the researches of other topics in logic since they can offer abstract representation for specific families of truth function operators. Apart of that, the conception of semi-expressiveness for abstract operators proposed in \textbf{Subsection \ref{S3}} might provide a new perspective to understand the underlying connection between representing a composite operator $R_{1}\circ\ldots\circ R_{n}$ and representing the operators $\lbrace R_{1},\ldots, R_{n}\rbrace$ respectively. It inspires us to address a new question: what is the sufficient and necessary condition for any composite operator $R_{1}\circ\ldots\circ R_{n}$ to be able to represent $\lbrace R_{1},\ldots, R_{n}\rbrace$? In the feature, we will continue to explore the relevant mechanisms of them.

Moreover, we suggest that our theory has the potential to be applied in the data science. Recently, it has received increasing attention to abstract and analyse logic structures from various data sets \cite{chen2012fuzzy,levy2000logic,kacprzyk2001linguistic}. There might be lots of potential truth function operators detected in the experiments. To verify their functional completeness (this determines how many logic relations the data contains), a purely constructive verification proof can be inefficient. Based on our method, their functional completeness can be verified in a more unified process, which makes it possible to design the computer-aided proof.

In a word, we suggest that the theory proposed in our work has the potential to be further explored.






\bibliographystyle{elsarticle-num}
\bibliography{Ref}

\begin{thebibliography}{10}
\expandafter\ifx\csname url\endcsname\relax
  \def\url#1{\texttt{#1}}\fi
\expandafter\ifx\csname urlprefix\endcsname\relax\def\urlprefix{URL }\fi
\expandafter\ifx\csname href\endcsname\relax
  \def\href#1#2{#2} \def\path#1{#1}\fi

\bibitem{shoenfield2018mathematical}
J.~R. Shoenfield, Mathematical logic, CRC Press, 2018.

\bibitem{radzki2017axiom}
M.~M. Radzki, On axiom systems of s{\l}upecki for the functionally complete
  three-valued logic, Axiomathes 27~(4) (2017) 403--415.

\bibitem{belnap1970every}
N.~Belnap, S.~McCALL, et~al., Every functionally complete m-valued logic has a
  post-complete axiomatization., Notre Dame Journal of Formal Logic 11~(1)
  (1970) 106.

\bibitem{venema1990expressiveness}
Y.~Venema, et~al., Expressiveness and completeness of an interval tense logic.,
  Notre Dame journal of formal logic 31~(4) (1990) 529--547.

\bibitem{wang2016functionally}
Z.-R. Wang, Y.-T. Su, Y.~Li, Y.-X. Zhou, T.-J. Chu, K.-C. Chang, T.-C. Chang,
  T.-M. Tsai, S.~M. Sze, X.-S. Miao, Functionally complete boolean logic in
  1t1r resistive random access memory, IEEE Electron Device Letters 38~(2)
  (2016) 179--182.

\bibitem{lehtonen2009stateful}
E.~Lehtonen, M.~Laiho, Stateful implication logic with memristors, in: 2009
  IEEE/ACM International Symposium on Nanoscale Architectures, IEEE, 2009, pp.
  33--36.

\bibitem{varshavsky2004functionally}
V.~Varshavsky, V.~Marakhovsky, I.~Levin, N.~Kravchenko, Functionally complete
  element for fuzzy control hardware implementation, in: The 2004 47th Midwest
  Symposium on Circuits and Systems, 2004. MWSCAS'04., Vol.~3, IEEE, 2004, pp.
  iii--263.

\bibitem{williamson1997resonant}
W.~Williamson~III, B.~K. Gilbert, Resonant tunneling diode structures for
  functionally complete low power logic, uS Patent 5,698,997 (Dec.~16 1997).

\bibitem{yetter1993functionally}
J.~D. Yetter, Functionally complete family of self-timed dynamic logic
  circuits, uS Patent 5,208,490 (May~4 1993).

\bibitem{tokmen1978functionally}
V.~Tokmen, A functionally-complete ternary system, Electronics Letters 14~(3)
  (1978) 69--71.

\bibitem{graham1967eta}
R.~L. Graham, On $\eta$-valued functionally complete truth functions, The
  Journal of Symbolic Logic 32~(2) (1967) 190--195.

\bibitem{massey1967binary}
G.~J. Massey, Binary connectives functionally complete by themselves in s5
  modal logic, The Journal of Symbolic Logic 32~(1) (1967) 91--92.

\bibitem{pinkava1978class}
V.~Pinkava, On a class of functionally complete multi-valued logical calculi,
  Studia Logica: An International Journal for Symbolic Logic 37~(2) (1978)
  205--212.

\bibitem{massey1970binary}
G.~J. Massey, et~al., Binary closure-algebraic operations that are functionally
  complete., Notre Dame Journal of Formal Logic 11~(3) (1970) 340--342.

\bibitem{maksimovic2006simple}
P.~Maksimovi{\'c}, P.~Jani{\v{c}}i{\'c}, Simple characterization of
  functionally complete one-element sets of propositional connectives,
  Mathematical Logic Quarterly 52~(5) (2006) 498--504.

\bibitem{brady2000peirce}
G.~Brady, From Peirce to Skolem: a neglected chapter in the history of logic,
  Elsevier, 2000.

\bibitem{chen2012fuzzy}
G.~Chen, Fuzzy logic in data modeling: semantics, constraints, and database
  design, Vol.~15, Springer Science \& Business Media, 2012.

\bibitem{levy2000logic}
A.~Y. Levy, Logic-based techniques in data integration, in: Logic-based
  artificial intelligence, Springer, 2000, pp. 575--595.

\bibitem{kacprzyk2001linguistic}
J.~Kacprzyk, R.~R. Yager, Linguistic summaries of data using fuzzy logic,
  International Journal of General System 30~(2) (2001) 133--154.

\end{thebibliography}
\end{document}